\documentclass{amsart}
\usepackage{lipsum}
\usepackage{amsfonts}
\usepackage{graphicx}
\usepackage{epstopdf}
\usepackage{mathtools}
\usepackage{float}
\usepackage{amsmath}
\usepackage{amssymb}
\usepackage{enumitem}
\usepackage{tikz}
\usepackage{booktabs}
\usepackage{makecell}
\usepackage{subcaption}
\usepackage{multirow}
\usepackage[
    backend=biber,
    style=alphabetic,
    sorting=ynt,
    maxnames=99,
    giveninits=true
]{biblatex}

\theoremstyle{plain}
\newtheorem{theorem}{Theorem}[section]
\newtheorem{lemma}[theorem]{Lemma}

\newtheorem{corollary}[theorem]{Corollary}
\newtheorem{proposition}[theorem]{Proposition}

\theoremstyle{definition}
\newtheorem{definition}[theorem]{Definition}

\theoremstyle{remark}
\newtheorem{remark}[theorem]{Remark}
\newtheorem{example}[theorem]{Example}

\title{Minimality of Tree Tensor Network Ranks}
\author{Jana Jovcheva}
\address{UCLouvain, INMA, ICTEAM, 1348 Louvain-la-Neuve, Belgium}
\email{jana.jovcheva@uclouvain.be}

\author{Tim Seynnaeve}
\address{Institute of Mathematics of the Polish Academy of Sciences, ul. {\'S}niadeckich 8, 00-656 Warszawa, Poland}
\email{tseynnaeve@impan.pl}

\author{Nick Vannieuwenhoven}
\address{KU Leuven, Department of Computer Science, Celestijnenlaan 200A, 3001 Leuven, Belgium}
\email{nick.vannieuwenhoven@kuleuven.be}

\begin{document}


\begin{abstract}
For a given tree tensor network $G$, we call a tuple of bond dimensions minimal if there exists a tensor $T$ that can be represented by this network but not on the same tree topology with strictly smaller bond dimensions. We establish necessary and sufficient conditions on the bond dimensions of a tree tensor network to be minimal, generalising a characterisation of Carlini and Kleppe about existence of tensors with a given multilinear rank. We also show that in a minimal tree tensor network, the non-minimal tensors form a Zariski closed subset, so minimality is a generic property in this sense.
\end{abstract}

\maketitle

\section{Introduction}
Tensor decompositions are powerful tools in modern machine learning, particularly deep learning, where they enable compact representations of high-dimensional data and model parameters. Structured decompositions like tensor trains and the Tucker format reduce the number of learnable parameters while maintaining expressive power \cite{Sidiropoulos_2017, stoundenmire}. These methods factor large tensors into networks of smaller core tensors connected by intermediate dimensions, known as \textit{bond dimensions}, which act as hyperparameters constraining the class of representable tensors. If the bond dimensions $\mathbf r$ are not minimal, the parameterisation is redundant, leading to wasted storage and computation. In practical algorithms such as rank-based truncation procedures \cite{grasedyckIntroductionHierarchicalRank2011}, density matrix renormalisation group (DMRG) methods \cite{SCHOLLWOCK201196}, or Riemannian optimisation on fixed-rank manifolds \cite{holtzManifoldsTensorsFixed2012, novikov2022}, they determine the size of local optimisation problems and the cost of intermediate contractions, so non-minimal choices increase computational complexity without enlarging the represented tensor class.

A natural question follows: \textit{what are the minimal bond dimensions required to exactly represent a given tensor class without redundancy?} The tuple of such minimal bond dimensions is known as the \textit{tree tensor network rank}.

\begin{example}\label{ex:overparameterisation}
As a familiar, basic example of redundancy or \textit{overparameterisation}, consider a matrix $M\in \mathbb R^{m\times n}$ that can be factorised as $M=XYZ$, where $X\in\mathbb R^{m\times r_1}$, $Y\in\mathbb R^{r_1\times r_2}$, $Z\in\mathbb R^{r_2\times n}$, and $r_1, r_2 \ll m, n$.  This matrix corresponds to the line tensor network depicted in Figure \ref{fig:line_tn}. If $M$ can be exactly represented by replacing $r_1, r_2$ with $q_1<r_1$ or $q_2<r_2$, then the original factorisation is overparameterised. 

\begin{figure}[tb]
    \centering
    \includegraphics[width=0.5\linewidth]{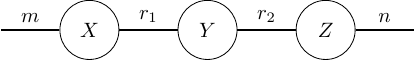}
    \caption{Line tensor network with bond dimensions $\mathbf{r}=(r_1, r_2)$ representing a matrix $M\in\mathbb R^{m\times n}$ factorised as $M=XYZ$.}
     \label{fig:line_tn}
\end{figure}

The row and column rank of a matrix are always equal, so a line tensor network model is not minimal if $r_1\ne r_2$, and the same set of matrices can be represented with $r_1'=r_2'=\min\{r_1,r_2\}$. In other words: $(r_1,r_2)$ is a tensor network rank if and only if $r_1=r_2$.

\end{example}

While different factorisations are possible, e.g., $M=A B$ with $A\in\mathbb R^{m\times r}$, $B\in\mathbb R^{ r\times n}$, and $r\ll m, n$, this work focuses on optimising bond dimensions within a \textit{fixed} tensor network structure rather than comparing factorisation schemes. Furthermore, our goal is not to approximate or analyse a specific tensor, but to \textit{characterise when a tensor network with fixed topology and bond dimensions defines a class of tensors for which no strictly smaller bond dimensions suffice}.

\begin{example}\label{ex:nonmin_vs_min} Consider, for example, the tree in Figure \ref{fig:min_vs_nonmin_tree} with two different tuples of bond dimensions.
\begin{figure}[tb]
    \centering
    \includegraphics[width=0.6\linewidth]{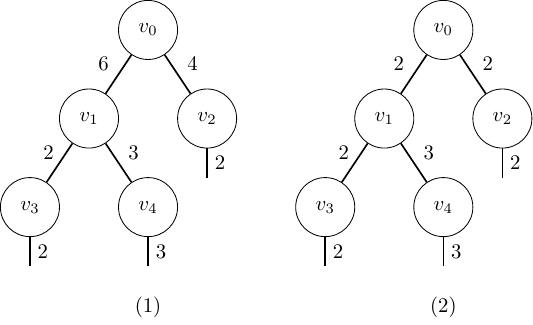}
    \caption{Does either of these labelings with bond dimensions correspond to a tree tensor network rank?}
    \label{fig:min_vs_nonmin_tree}
\end{figure}
Our main result, Theorem \ref{theorem:minimality}, implies that all tensors representable with (1) are also representable with (2). Since the latter bond dimensions are componentwise less than the former, the labeling of the left tree is not minimal. Therefore, representing tensors with (1) includes redundancy. Theorem \ref{theorem:minimality} will show that the left panel is not minimal because the dimensions on bonds $(v_0, v_1)$ and $(v_0, v_2)$ exceed the product of the dimensions on adjacent bonds, while the right panel is minimal because all bond dimensions are admissible (Definition \ref{def:admissible}).
\end{example}

\subsection{Related work}\label{sec:lit_review}

Tensor networks originated in quantum physics as a means to represent many-body systems, with entanglement encoded in the underlying graph \cite{Verstraete01032008, ShiDuanVidal2006, evenblyTensorNetworkStates2011, Orus2014}. In this context, Shi, Duan, and Vidal established a canonical form for tree tensor networks based on singular value decomposition (SVD) in which every edge is associated with a Schmidt decomposition across the bipartition of the tree obtained by cutting that edge. In this form, the minimal bond dimension on an edge equals the (Schmidt) rank across the cut \cite{ShiDuanVidal2006}. 

In applied mathematics, there emerged a differential-geometric description of matrix product states/tensor trains and tree tensor networks/hierarchical Tucker tensors, which form smooth embedded manifolds once the ranks are fixed and certain full-rank (or ``minimal'') conditions are imposed \cite{holtzManifoldsTensorsFixed2012, USCHMAJEW2013133, Haegeman2014GeometryMPS}. These manifold constructions describe the geometry of a fixed-rank model class modulo internal gauge symmetries. Crucially, all of these works assume that the tuple of bond dimensions is minimal, i.e., a tree tensor network rank. However, this literature does not provide explicit criteria for verifying minimality. In Figure \ref{fig:min_vs_nonmin_tree}, both sets of bond dimensions may appear valid, but (1) is not minimal. Since the manifold dimension is determined by the minimal bond dimensions, a non-minimal tuple merely yields a redundant parameterisation. Verifying minimality is therefore essential for the differential-geometric picture.

In applied algebraic geometry, sets of tensors representable with bounded bond dimensions are studied as algebraic varieties defined by polynomial equations. In this setting, a fixed tensor network topology together with a bond dimension tuple determines a parameterised family of tensors whose Zariski closure forms a tensor network variety. Prior work has investigated structural properties of these varieties, including dimension and geometric features \cite{landsbergGeometryTensorNetwork2012, buczynskaHackbuschConjectureTensor2015, limandye, bernardiDimensionTensorNetwork2023}. In particular, in \cite{bernardiDimensionTensorNetwork2023} it is shown that overabundant bond dimensions can be reduced without affecting the dimension of the associated tensor network variety. Our results fit naturally into this framework: minimality becomes a statement about when generic tensors in the associated tensor network variety do not admit a representation with any strictly smaller bond dimensions. 

Existing works have characterised the tree tensor network rank in terms of flattenings of bipartitions of the tree induced by edge cuts. Hackbusch and Kühn introduced the hierarchical Tucker format together with a recursive SVD-based truncation procedure that enables controlled bond dimension reduction \cite{hackbuschNewSchemeTensor2009a}. Grasedyck and Hackbusch subsequently formalised hierarchical and TT ranks via dimension trees and related them to matrix ranks of edge-induced flattenings, establishing structural relations between the two formats \cite{grasedyckIntroductionHierarchicalRank2011}. Bachmayr, Schneider, and Uschmajew characterised representability in hierarchical tensor formats through these flattening ranks, a viewpoint that underlies the hierarchical SVD method \cite[Theorem 3.3]{bachmayrTensorNetworksHierarchical2016a}. Ye and Lim showed that the tree tensor network rank of a tensor is unique and can be computed as the tuple of ranks of global edge-cut flattenings \cite[Theorems 8.3 and 8.8]{limandye}.

Most closely related to this article are the works by Carlini and Kleppe \cite{carlinikleppe} and Buczy\'nska, Buczy\'nski, and Micha{\l}ek \cite{buczynskaHackbuschConjectureTensor2015}. The former determined explicit inequalities ensuring that a tuple of bond dimensions is minimal for star graphs. The latter determined necessary conditions ensuring that a tuple of bond dimensions is minimal for binary trees, formulating them recursively via an auxiliary function, but did not prove sufficiency.

\subsection{Tensor networks}\label{sec:prelim}
To enable a precise statement of our contributions, we now recall the formal definition of tensor networks for arbitrary fields from the work by Ye and Lim \cite{limandye}. In this paper, we are specifically interested in \textit{tree tensor networks}, whose underlying \textit{contraction graph} is a connected acyclic graph. Let $G = (V, E)$ be a tree, where $V=\{1,\dots, d\}$ are vertices and $E \subseteq \{(i, j) : i, j \in V, i \ne j\}$ are directed edges. The ordering $(i,j)$ indicates that the edge is from vertex $i$ to vertex $j$. Each (undirected) edge has a corresponding nonzero bond dimension $r_{ij} \in \mathbb N$. Each vertex $i\in V$ has an associated ``physical'' vector space $\mathbb V_i$, which is of finite dimension over a field $\Bbbk$. 

\begin{figure}[tb]
    \centering
    \includegraphics[width=0.42\linewidth]{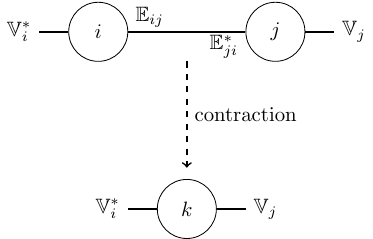}
    \caption{Contraction of two tensor vertices $i$ and $j$ along shared edge $(i, j)$, pairing the vector space $\mathbb E_{ij}$ with its dual $\mathbb E_{ji}^*$, resulting in a new node $k$ with structure inherited from $\mathbb V_i^*$ and $\mathbb V_j$. Note that the dangling edges are not edges in $E$, and merely provide a visualisation of the vector space associated with each vertex.}
    \label{fig:contraction}
\end{figure}

\begin{remark}
	The physical spaces $\mathbb V_i$ will be visualised with ``dangling'' edges, as in Figure \ref{fig:contraction}, though they are distinct from the bond edges in $E$ which encode the contraction structure of the network. In the literature it is often allowed for a vertex to have more than one physical output associated with it, or none at all. Restricting to exactly one physical space at each vertex does not limit generality: at vertices with several physical outputs, we define $\mathbb V_i$ to be the tensor product of the constituent vector spaces. At vertices with no physical output, we can take $\mathbb V_i=\Bbbk$. Since $\Bbbk\otimes W\cong W$, the presence of such 1-dimensional physical spaces has no impact on the network. 
\end{remark}

A tensor network will be used to represent a tensor in the tensor product space $\mathbb V_1\otimes\dots\otimes\mathbb V_d$ implicitly as the result of \textit{contracting} the network as specified by the contraction graph. To this end, each edge $(i,j)$ has associated ``contraction'' vector and covector spaces $\mathbb E_{ij}$ and $\mathbb E_{ji}^*$, respectively, of dimension equal to the bond dimension $r_{ij}$. Here $\mathbb E_{ji}^*$ denotes the dual space of $\mathbb E_{ij}$, i.e., the space of all linear maps $\phi^* : \mathbb E_{ij} \to \Bbbk$. Each map $\phi^* \in \mathbb E_{ji}^*$ evaluates the vector $\psi \in \mathbb E_{ij}$ via $ \phi^*(\psi)$, where $\phi^*(\psi)$ denotes the action of $\phi^*$ on $\psi$. This pairing underlies the contraction operations used in the network.

\begin{example}\label{ex:contraction}
To illustrate the contraction process, consider the two connected vertices $i$ and $j$ associated with vector spaces $\mathbb V_i^*$ and $\mathbb V_j$, respectively, in Figure \ref{fig:contraction}. The top diagram represents a contraction along the shared edge $(i,j)$, which can be interpreted as matrix multiplication: the left node encodes a tensor $A \in \mathbb V_i^* \otimes \mathbb E_{ij} = \mathrm{Hom}(\mathbb V_i, \mathbb E_{ij})$, and the right node encodes $B \in \mathbb E_{ji}^* \otimes \mathbb V_j = \mathrm{Hom}(\mathbb E_{ij}, \mathbb V_j)$, where $\mathrm{Hom}(\mathbb V, \mathbb W)$ is the space of linear maps from $\mathbb V$ to $\mathbb W$. Contracting along the shared edge corresponds to composing these two maps: $C = BA \in \mathrm{Hom}(\mathbb V_i, \mathbb V_j)$. The resulting node $k$ in the bottom diagram inherits the physical vector spaces $\mathbb V_i^*$ and $\mathbb V_j$, and the corresponding tensor $C$ encodes the composition of $A$ and $B$, i.e., the result of contracting $\mathbb E_{ij}$ with $\mathbb E_{ji}^*$.
\end{example}

The above example generalises to arbitrary tensor contractions in tree tensor networks, as explained next. Recall that a tensor $v\in\mathbb V_1\otimes\dots\otimes\mathbb V_d$ is \textit{elementary} if it can be written in the form $v_1\otimes\dots\otimes v_d$ for $v_i\in\mathbb V_i$. Any tensor can be expressed as a finite sum of elementary tensors, $T=\sum_{j=1}^s v_{1}^{(j)}\otimes\dots\otimes v_{d}^{(j)}$ \cite[Section 3]{Greub1978}. For each vertex $i\in V$, let $\text{in}(i)$ and $\text{out}(i)$ denote the set of incoming and outgoing edges incident to $i$. Let the \textit{inspaces} and \textit{outspaces} be 
\[
    \mathbb I_i^*=\bigotimes _{j\in\text{in}(i)}\mathbb E^*_{ji} \quad \text{and} \quad \mathbb O_i=\bigotimes_{j\in\text{out}(i)}\mathbb E_{ij},
\]
respectively. Then, a tree $G$ represents a tensor $T$ via the \textit{contraction map}
\begin{equation}\label{eq:contraction_map}
g: \bigotimes_{i=1}^d
\left(\mathbb I_i^*
\otimes \mathbb V_i \otimes \mathbb O_i \right) 
\rightarrow \mathbb V_1\otimes\dots\otimes\mathbb V_d,
\end{equation}
which contracts each pair of dual vector spaces $\mathbb E_{ji}^* \otimes \mathbb E_{ij}$ along the directed edge $(i, j) \in E$. The map $g$ is linear and defined by its action on elementary tensors: the contraction
\[
g\left(\bigotimes_{i=1}^d(\phi^*_{i} \otimes v_i \otimes \psi_{i})\right) \coloneq \prod_{(i, j)\in E} \phi^*_{ij}(\psi_{ij})\cdot \bigotimes_{i=1}^d v_i,
\]
where $v_i \in \mathbb V_i$, $\phi^*_{i}=\bigotimes_{j\in\mathrm{in}(i)}\phi_{ij}^*, \phi^*_{ij} \in \mathbb E^*_{ji}$, and $\psi_{i}=\bigotimes_{j\in\mathrm{out}(i)}\psi_{ij}, \psi_{ij} \in \mathbb E_{ij}$, acts by evaluating the pairings $\phi^*_{ij}(\psi_{ij})$ for each edge $(i, j)$. 

The set of all tensors in $\mathbb V_1\otimes\dots\otimes \mathbb V_d$ representable by the contraction graph $G$ and the tuple of bond dimensions $\mathbf{r}=(r_{ij}=\dim\mathbb E_{ij})_{(i,j)\in E}$ is denoted by $\mathrm{TN}(G, \mathbf r)$.  That is, $\mathrm{TN}(G,\mathbf{r})$ contains the tensors that can be written as a contraction $g(T_1 \otimes \dots \otimes T_d)$,
where each \textit{local tensor} $T_i \in \mathbb I_i^* \otimes \mathbb V_i \otimes \mathbb O_i$, i.e.,
\[
    \mathrm{TN}(G, \mathbf{r})=\{g(T_1\otimes\dots\otimes T_d) \mid T_i\in \mathbb I_i^*\otimes \mathbb V_i\otimes \mathbb O_i\}.
\]

\subsection{Contributions}
This work considers the minimality of bond dimensions. This concept is defined next. For two bond dimension tuples $\mathbf{r} = (r_{ij})_{(i,j) \in E}$ and $\mathbf{s}= (s_{ij})_{(i,j) \in E}$, we write $\mathbf{s} \leq \mathbf{r}$ to mean componentwise inequality $s_{ij} \leq r_{ij}  \ \forall (i,j) \in E$. 

\begin{definition}[Minimal bond dimensions]\label{definition:minimal_bond_dims}
    Fix a tree $G$ and a tuple of bond dimensions $\mathbf{r}$.
    \begin{itemize}
    \item Given $T \in \mathrm{TN}(G,\mathbf{r})$, we say that $\mathbf{r}$ is \emph{minimal for $T$}, or a \emph{tensor network rank of $T$}, if  $T \in \mathrm{TN}(G, \mathbf{r})$, but $T \notin \mathrm{TN}(G, \mathbf{s})$ for any $\mathbf{s} \leq \mathbf{r}$ with $\mathbf{s} \neq \mathbf{r}$.  
    \item We say that $\mathbf{r}$ is  \emph{minimal for $G$}, or a \emph{tensor network rank for $G$} if there exists a $T \in \mathrm{TN}(G,\mathbf{r})$ such that $\mathbf r$ is minimal for $T$.
    \end{itemize}
\end{definition}

Consider the set
\[
    \mathrm{TN}^\circ(G, \mathbf r)=\mathrm{TN}(G, \mathbf r) \setminus \bigcup_{\begin{smallmatrix} \mathbf s\le\mathbf r \\  \mathbf s\ne\mathbf r \end{smallmatrix}} \mathrm{TN}(G, \mathbf s).    
\]
By definition, $\mathrm{TN}^{\circ}(G, \mathbf{r})$ consists of all tensors $T$ such that $\mathbf r$ is minimal for $T$. So $\mathbf r$ is minimal for $G$ if and only if $\mathrm{TN}^{\circ}(G, \mathbf{r})$ is nonempty. Note that $\mathrm{TN}^\circ(G, \mathbf{r})$ can be the empty set, as was shown in Example \ref{ex:overparameterisation} when $r_1\ne r_2$. 

The main questions we resolve in this paper are the following:

\begin{enumerate}[label=Q\arabic*]
    \item \textbf{Minimality:} When is $\mathrm{TN}^{\circ}(G,\mathbf{r})$ nonempty, i.e., when is $\mathbf r$ a tensor network rank for $G$?

    \item \textbf{Genericity:} If $\mathrm{TN}^\circ(G, \mathbf r)\ne \emptyset$, is it a Zariski open and dense subset of $\mathrm{TN}(G, \mathbf r)$? That is, does minimality imply that $\mathrm{TN}^\circ(G, \mathbf r)$ contains almost all tensors of $\mathrm{TN}(G, \mathbf r)$?
\end{enumerate}

Much of the literature, e.g. \cite{holtzManifoldsTensorsFixed2012, USCHMAJEW2013133, Haegeman2014GeometryMPS}, proves results under the assumption that $\mathbf r$ is a tree tensor network rank. The main result of this work, Theorem \ref{theorem:minimality}, provides a necessary and sufficient condition on $(G, \mathbf r)$ for $\mathbf r$ to be a tree tensor network rank for $G$, generalising the main result of Carlini and Kleppe \cite{carlinikleppe} to arbitrary tree tensor networks. When this condition holds, we show in Theorem \ref{theorem:genericity} below that the set of tensors $\mathrm{TN}^\circ(G, \mathbf r)$ even forms a \textit{Zariski open and dense} subset of $ \mathrm{TN}(G, \mathbf{r}) $. 

Our arguments are purely algebraic, do not rely on orthogonality or singular value decompositions, and apply over arbitrary infinite ground fields. In this way, our work contributes a precise local-to-global minimality characterisation and genericity statement for tree tensor network varieties, with easy-to-verify conditions that involve only the numerical values in a tuple of bond dimensions.

\subsection{Outline}
The remainder of this paper is structured as follows. We begin in Section \ref{sec:star_graphs} by recalling Carlini and Kleppe's characterisation of the star graph case, which corresponds to Tucker decompositions, and setting the algebraic-geometric foundation that underpins our generalisation to arbitrary trees. Thereafter, in Section \ref{sec:main}, we extend this discussion to arbitrary trees and state the necessary and sufficient conditions for minimality on these trees. We then establish the genericity of these minimal models within the ambient tensor network variety, and conclude by discussing practical reduction of non-minimal tree tensor networks.

\section{Carlini and Kleppe's characterisation for star graphs}\label{sec:star_graphs}
\begin{figure}[tb]
    \centering
    \includegraphics[width=0.35\linewidth]{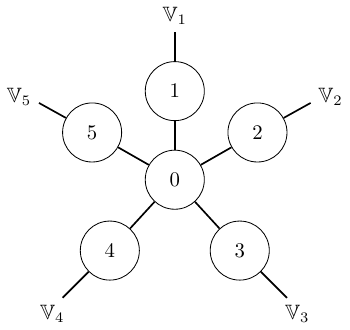}
    \caption{Star graph for a 5th order tensor. }
    \label{fig:star-graph}
\end{figure}

A \textit{Tucker graph} (or star graph) is a tensor network with a single central node connected to $d$ leaf nodes, where $d$ is the order of the tensor (see Figure \ref{fig:star-graph}). Formally, a Tucker graph $G_{\text{Tucker}}$ has vertices $V = \{0, 1, \dots, d\}$, where $0$ is the central node, and edges $E = \{(0, i) : i = 1, \dots, d\}$ \cite{limandye}. We require the space associated with vertex $0$ to be trivial: $\mathbb{V}_0=\Bbbk$. Then, the resulting tensor network represents the Tucker decomposition of an order-$d$ tensor. 

Carlini and Kleppe \cite{carlinikleppe} provided a complete characterisation of minimality (Q1) and genericity (Q2) for Tucker graphs by establishing necessary and sufficient conditions on the bond dimensions. We briefly recall this background and use this discussion to introduce the mathematical techniques and notation we need to prove the corresponding result for general tree tensor networks in the next section.

Recall that a \textit{flattening}, or \textit{matrix unfolding}, is the linear mapping 
\begin{equation}\label{eq:flattening}
    \bullet_{(S)}: \mathbb V_1\otimes\dots\otimes \mathbb V_d\rightarrow \left(\bigotimes_{j\in S}\mathbb V_j\right)\otimes\left(\bigotimes_{k\in S^c}\mathbb V_k\right)
\end{equation}
associated with the partition of a tensor into two disjoint subsets, $S\cup S^c=\{1,\dots,d\}$, where $S\subseteq \{1,\dots,d\}$ and $S^c$ is its complement. This mapping, written as $T_{(S)}$, sends an elementary tensor $T=v_1\otimes\dots\otimes v_d $ to the matrix
\[
    T_{(S)}=\left(\bigotimes_{j\in S} v_j\right)\otimes\left(\bigotimes_{k\in S^c} v_k\right)\in\mathrm{Hom}\left(\bigotimes_{k\in S^c}\mathbb V_k^*,\bigotimes_{j\in S}\mathbb V_j\right).
\]
By the universal property of the tensor product \cite[Section 1.4]{Greub1978}, this fully defines $\bullet_{(S)}$. Each choice of partition $S$ defines a distinct unfolding. If $S=\{m\}$ and $S^c = \{1,\dots,d\} \setminus \{m\}$, the resulting matrix is the \textit{mode-$m$ unfolding} of $T$, denoted more succinctly by $T_{(m)}$. The rank of this matrix is known as the \textit{flattening rank} along mode $m$ of the tensor. In the context of a tree tensor network $G=(V, E)$, we will frequently consider bipartitions induced by edge cuts. For an edge $(i,j)\in E$, removing $(i,j)$ disconnects $G$ into two subtrees, resulting in a partition $S_{ij}\cup S_{ij}^c=\{1,\dots,d\}$ of the tensor modes. We denote the corresponding flattening by $T_{(i,j)}=T_{(S_{ij})}$. Note that while the flattening map formally depends on a sequence of modes, we describe the bipartition using sets for notational convenience. This is justified because the flattening rank is invariant under permutations of the tensor modes. 

The \textit{multilinear rank} of a tensor $T$ is the tuple $\boldsymbol{\mu}=(\mu_1,\dots, \mu_d)$, where $\mu_i = \mathrm{rank}(T_{(i)})$ is the rank of the mode-$i$ unfolding as defined by \eqref{eq:flattening}. It is well known that for our Tucker graph, the set $\mathrm{TN}(G_\text{Tucker}, \boldsymbol{r})$ consists precisely of the tensors with multilinear rank bounded by $\boldsymbol{r} = (r_{01},\ldots,r_{0d})$, and moreover it is an \emph{irreducible affine algebraic variety} \cite{carlinikleppe, Landsberg2012}. 

Recall from \cite{idealsvarieties} that an \textit{affine variety} $\mathcal V$ in the affine space $\Bbbk ^N$ is a set of solutions to a system of polynomial equations $\{f_1,\dots, f_n\}$ in $\Bbbk$, i.e.,
\[
\mathcal V(f_1,\dots,f_n)=\{(v_1,\dots,v_N)\in \Bbbk^N: f_i(v_1,\dots, v_N)=0,\quad \forall 1\le i\le n\}.
\]
Affine varieties form the closed sets of a topology on $\Bbbk ^N$, called the \emph{Zariski topology}. Complements of affine varieties are called \emph{Zariski open}.
By fixing bases in each $\mathbb V_i$, we identify the tensor product $\mathbb V_1 \otimes \dots \otimes \mathbb V_d$ with the affine space $\Bbbk^{n_1 \cdots n_d}$, where $n_i = \dim \mathbb V_i$. In this basis, each tensor is represented by its coordinate entries $T_{i_1,\dots,i _d}$, which are interpreted as variables in the polynomial ring $\Bbbk[T_{i_1,\dots, i_d}]$, where $1\le i_j\le n_j$. For a tensor $T\in\mathbb V_1\otimes\dots\otimes\mathbb V_d$, the condition that its unfolding along mode $i$ has rank at most $\mu_i$ is equivalent to the vanishing of all $(\mu_i+1)\times(\mu_i+1)$ minors of the corresponding matrix flattening. These minors are polynomial functions of the entries of $T$, so the set of all tensors with multilinear rank at most $(\mu_1,\dots,\mu_d)$ is exactly the set of tensors for which these polynomials vanish. In other words, this set is defined by a system of polynomial equations over $\Bbbk$, and thus forms an algebraic variety. 

An affine variety is called \emph{irreducible} if it cannot be written as a union of two proper subvarieties. Over an infinite field, the image of an affine space under a polynomial map is irreducible \cite[Proposition 4.5.5]{idealsvarieties}. Hence we obtain the following lemma.  

\begin{lemma}\label{lem:irreducibility}
    The variety $\mathrm{TN}(G, \boldsymbol{r})$ is irreducible for all bond dimensions $\mathbf{r}$ and every contraction graph $G$ that is a tree.
\end{lemma}

This is shown by Lim and Ye \cite[Proposition 2.5]{limandye}, where $\mathrm{TN}(G, \mathbf{r})$ is realised as the image of the polynomial contraction map from the product of vertex and edge tensor spaces. 

The main result of Carlini and Kleppe, \cite[Theorem 7]{carlinikleppe}, states the set of star graphs $\mathrm{TN}^\circ(G_\text{Tucker}, \mathbf r)$ is nonempty if and only if the bond dimensions satisfy
\begin{equation} \label{eq:TuckerMinimal}
r_{0i} \le \prod_{j \ne i} r_{0j} \quad \text{ and } \quad r_{0i} \le \dim \mathbb V_i \quad \forall \ i = 1,\dots,d .
\end{equation}

Since $\mathrm{TN}(G, \mathbf r)$ is irreducible, any nonempty Zariski open subset is Zariski dense. Therefore, if $\mathrm{TN}^\circ(G, \mathbf r)$ is nonempty, it will be Zariski open and dense in $\mathrm{TN}(G, \mathbf r)$, meaning that a generic tensor in $\mathrm{TN}(G, \mathbf r)$ will satisfy these conditions. To summarise, the characterisation by Carlini and Kleppe \cite{carlinikleppe} provides definitive answers to both of our questions for Tucker graphs:

\begin{enumerate}[label=A\arabic*]
    \item \textbf{Minimality:} $\mathrm{TN}^\circ(G_{\text{Tucker}}, \mathbf{r}) \neq \emptyset$ if and only if $\mathbf{r}$ satisfies \eqref{eq:TuckerMinimal}.

    \item \textbf{Genericity:} If $\mathrm{TN}^\circ(G_{\text{Tucker}}, \mathbf{r}) \neq \emptyset$, then it is a Zariski open and dense subset of $\mathrm{TN}(G_{\text{Tucker}}, \mathbf{r})$. 
\end{enumerate}

\section{The main result on arbitrary tree graphs}\label{sec:main}
A star graph is a specific tree. We now ask more generally which constraints govern minimality across the full class of tree tensor networks. The goal of this section is to prove that necessary and sufficient conditions for minimality follow from enforcing the same type of bond dimension constraints on each local tensor $T_i$ in the graph.

We now fix some notation. We say that a subtree $G' = (V', E')$ of $G$ is \emph{rooted at $a \in V$} if every edge $(i,j) \in E'$ is oriented such that $j$ is closer to $a$. Then every vertex except for $a$ will have a single outspace:  $\mathbb O_{i}=\mathbb E_{i, \text{out}(i)}$. The dimensions of the spaces $\mathbb O_i$ correspond to the imposed bond dimensions $r_{i, \text{out}(i)}$. 
Let $\mathbb {\widetilde I}_i=\mathbb I_i^*\otimes \mathbb V_i$. For each $i\in V$, consider the flattening of the local tensor $T_i$ where $\mathbb {\widetilde I}_i$ is mapped to the rows of the matrix and $\mathbb O_i$ is mapped to the columns, which yield flattening matrices 
\begin{equation}\label{eq:injective_map}
    M_{i}:\mathbb O_i\rightarrow \widetilde{\mathbb I}_i.
\end{equation}
Provided that $\dim \mathbb {\widetilde I}_i\ge r_{\text{out}(i)}=\dim\mathbb O_i$, i.e., the number of columns in the flattening matrix is less than or equal to the number of rows, the matrix can have linearly independent columns, or equivalently, $M_i$ is the matrix of an injective linear map.  

The following lemma formalises how cutting an edge and considering locally maximal flattening ranks enforces bond dimension constraints within a connected subtree.

\begin{lemma}[Injectivity propagation]\label{lem:edge-rank}
    Let $T\in\mathbb V_1\otimes\dots\otimes\mathbb V_d$ be represented by a tree tensor network $G=(V, E)$.  Fix an edge $(a,p)\in E$, and let $G_a=(V_a,E_a)$ be the subtree obtained by removing $(a,p)$ and retaining the component containing $a$, viewing $(a,p)$ as a dangling edge with physical space $\mathbb E_{ap}$. Assume that $G_a$ is rooted at $a$, and that the flattening $M_i$ from \eqref{eq:injective_map} of every local tensor $T_i$ has rank exactly $r_{i,\mathrm{out}(i)}$. Denote by $\widehat M_a$ the $(\prod_{i \in V_a} \dim {\mathbb V}_i ) \times r_{ap}$-matrix obtained by contracting all internal edges in $G_a$. Then, $\mathrm{rank}(\widehat M_a) = r_{ap}$, i.e., $\widehat M_a$ represents an injective map.
\end{lemma}

\begin{proof}
We proceed by induction on the depth of the connected subtree $G_a$ rooted at $a$. The depth of a rooted tree is the number of edges on the longest simple path from the root $a$ to any leaf in $G_a$. The base case is the trivial subtree of depth zero consisting only of vertex $a$. In this case, $\widehat M_a=M_a$ is simply the flattening of the local tensor $T_a$ along the partition induced by edge $(a,p)$, which by assumption has rank $r_{ap}$. 

Now assume the lemma holds for all connected subtrees of depth strictly less than $d$. Let $G_a$ be a connected subtree of depth $d$ with root $a$, and let $\{c_1,\dots, c_n\}$ be its children. For each child $c_i$, let $G_{c_i}$ be the subtree rooted at $c_i$ which is of depth at most $d-1$. By the inductive hypothesis, $\mathrm{rank}\bigl(\widehat M_{c_i}\bigr)=r_{ac_i}$ for each $i$, and $\mathrm{rank}\left( M_a\right)=r_{ap}$. Recall that $g(\cdot)$ is the contraction map from \eqref{eq:contraction_map} relating the tree $G$ to the tensor $T$. Then, the flattening matrix resulting from a contraction of all internal edges in the subtree and flattening of the resulting tensor is
\[
    \widehat M_a = g\left(\left(\bigotimes_{i=1}^n \widehat M_{c_i}\right) \otimes  M_a\right).
\]
Since all local tensors are flattened, $g(\cdot)$ corresponds to the matrix multiplication
\[
        \widehat M_a = \left(\bigotimes_{i=1}^n \widehat M_{c_i}\right)  M_a.
\]
As $\widehat M_{c_i}$ has linearly independent columns by the inductive hypothesis, it defines an injective linear map. By \cite[equation 1.12]{Greub1978}, the tensor product of injective maps is injective, so $\bigotimes_{i=1}^n \widehat M_{c_i}$ is injective. Multiplying on the right by $M_a$, which is injective by assumption, preserves injectivity. Hence, $\mathrm{rank}(\widehat M_a)=\mathrm{rank}(M_a)=r_{ap}$. This concludes the proof. 
\end{proof}

\begin{remark}
Lemma \ref{lem:edge-rank} may also be interpreted through the normal form constructions for tree tensor networks obtained by suitable orthogonalisations of edge-induced flattenings (see, e.g., \cite[Section 11.3.2]{hackbusch2019tensor} and \cite{grasedyckIntroductionHierarchicalRank2011, Grasedyck2024}). In such gauges, each local tensor defines an isometric map toward its parent, and injectivity along an edge propagates automatically under contraction. The proof above emphasises that this result is purely algebraic rather than metric; it does not require an inner product.
\end{remark}

This establishes that the rank of the contracted tensor associated with a subtree is equal to the bond dimension along the edge connecting it to the rest of the network under an assumption on the flattening ranks of the local tensors. We now proceed to characterise minimal bond dimensions in tree tensor network representations. First, we need two definitions.

\begin{definition}[Admissible bond dimensions]\label{def:admissible}
	Let $G=(V,E)$ be a tree with vector spaces $\mathbb V_i$ at each vertex $i\in V$, and bond dimensions $r_{ij}=\dim\mathbb E_{ij}$ for each edge $(i,j)\in E$. The bond dimensions are called \emph{admissible} if, for every vertex $i$ and neighbour $j \in \mathrm{nb}(i)$, we have
	\begin{equation}\label{eq:admissible}
	r_{ij} \le \dim \mathbb V_i\prod_{k \in \mathrm{nb}(i)\setminus\{j\}} r_{ik}.
	\end{equation}
\end{definition}

\begin{definition}[Effective multilinear rank]\label{def:eff_ml_rank}
Let
\[
T_i \in \left(\!\bigotimes_{k\in\mathrm{in}(i)} \mathbb E_{ki}^*\right)
   \otimes \mathbb V_i \otimes
   \left(\!\bigotimes_{k\in\mathrm{out}(i)} \mathbb E_{ik}\right)
\]
be the local tensor at vertex $i$ of a tree tensor network $G=(V, E)$ representing a tensor $T$. For each neighbour $j \in \mathrm{nb}(i)$, let $T_i^{(j)}$ denote the matricisation of $T_i$ obtained by taking the bond space (i.e., $\mathbb E_{ij}$ or $\mathbb E_{ji}^*$) associated with edge $(i,j)$ as the row space, and all remaining incident bond spaces together with $\mathbb V_i$ as the column space. The \emph{effective multilinear rank} of $T_i$ is then the tuple
\[
\boldsymbol{\mu}^{\mathrm{eff}}(T_i)
= \bigl(\mathrm{rank} (T_i^{(j)}) \mid j\in\mathrm{nb}(i)\bigr).
\]
This tuple coincides with the multilinear rank of $T_i$ with respect to the decomposition into its incident bond spaces and the physical space $\mathbb V_i$, except that the component corresponding to the flattening separating $\mathbb V_i$ from the product of all bond spaces is omitted.
\end{definition}

 The following result is a slight generalisation of \cite[Theorem 7]{carlinikleppe}.
 
\begin{proposition}[Genericity of effective multilinear rank]\label{prop:fullEffectiveRank}
	If the admissibility inequalities \eqref{eq:admissible} are satisfied at a given vertex $i$, then the condition that a local tensor $T_i$ satisfies $\mathrm{rank} (T_i^{(j)}) = r_{ij}$ for every $j \in \mathrm{nb}(i)$ defines a nonempty Zariski open dense subset of the local space. In other words, a generic tensor will have full effective multilinear rank.
\end{proposition}

\begin{proof}
    For every $j$, the condition $\mathrm{rank}(T_i^{(j)}) = r_{ij}$ defines a Zariski open subset of the local space, given by the nonvanishing of suitable minors.  By \eqref{eq:admissible} and \cite[Theorem 7]{carlinikleppe}, there exists a tensor attaining all of these ranks simultaneously. Hence the intersection of the corresponding open sets is nonempty Zariski open (and thus dense).
\end{proof}

We will further require the following three results by Ye and Lim \cite{limandye}.

\begin{lemma}[Uniqueness of tree tensor network rank. {\cite[Theorem 8.3]{limandye}}] \label{lem:uniqueTTNrank}
    For a fixed tree $G$, the tree tensor network rank of a tensor $T\in\mathbb V_1\otimes\dots\otimes\mathbb V_d$ is unique.
\end{lemma}

\begin{lemma}[Tree tensor network rank characterisation. {\cite[Theorem 8.8]{limandye}}]\label{lem:TTNrank}
    The tree tensor network rank $\mathrm{rank}_G(T)=\mathbf r$ is given by
    \[
        r_{ij}=\mathrm{rank}(T_{(i,j)}),
    \]
    where $T_{(i,j)}$ is the flattening of $T$ along the edge $(i, j)$.
\end{lemma}

\begin{corollary}[$\mathrm{TN}(G, \mathbf{r})$ is an algebraic variety {\cite[Corollary 8.9]{limandye}}] \label{cor:TNvariety}
    $\mathrm{TN}(G, \mathbf{r})$ is an irreducible algebraic variety whose vanishing ideal is generated by all $(r_{ij}+1)\times(r_{ij}+1)$ minors of the flattening map \eqref{eq:flattening}.
\end{corollary}

While Ye and Lim work over the complex field, these statements remain valid for arbitrary infinite fields (see \cite{adams2024}), which justifies our general setting.

We are now ready to state and prove the first main result.

\begin{theorem}[Minimality]\label{theorem:minimality}
Let $T=g(T_1\otimes\dots\otimes T_d)$ be a tensor represented by a tree tensor network $G = (V, E)$ with bond dimensions $\mathbf{r} = \bigl(r_{ij} \mid (i, j)\in E\bigl)$. Then $\mathbf{r}$ is minimal for $T$ if and only if, for every vertex $i \in V$, the local tensor $T_i$ has effective multilinear rank equal to the bond dimensions $\bigl(r_{ij} \mid j\in \mathrm{nb}(i)\bigl)$ along its incident edges.
\end{theorem}

\begin{proof}[Proof] \textbf{Necessity.} Suppose there exists an internal vertex $i \in V$ with incident edge $(i, j) \in E$ such that $\mu_{ij}\coloneq\mathrm{rank}\bigl(T_i^{(j)}\bigr)<r_{ij}$, using the notation from Definition \ref{def:eff_ml_rank}. From Section \ref{sec:star_graphs}, $T_i$ can be represented as a star graph. Consider the construction illustrated in Figure \ref{fig:tucker-trick}. The mode-$j$ rank of ${T}_i$ is $\mu_{ij} < r_{ij}$, so the factor matrix $A_i^{(j)}$ has fewer columns than the bond dimension. 

\begin{figure}[tb]
    \centering
    \includegraphics[width=0.8\linewidth]{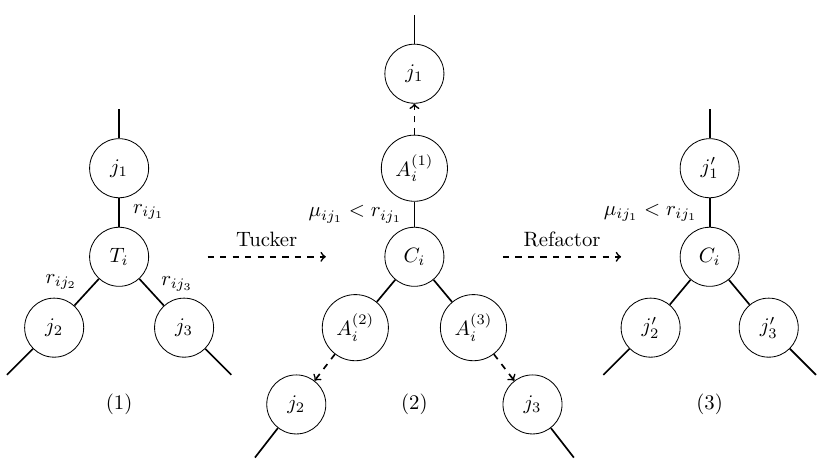}
    \caption{Illustration of the refactoring procedure. (1) Local tensor at vertex $i$ connected to neighbouring vertices $j_1, j_2, j_3$ with bond dimensions $r_{ij_1}, r_{ij_2}, r_{ij_3}$. (2) Tucker decomposition of ${T}_i$ into core tensor $C_i$ and factor matrices $A_i^{(j)}$. (3) Refactored network with reduced bond dimension $\mu_{ij_1}<r_{ij_1}$ achieved by absorbing factor matrices into vertices adjacent to $i$.}
    \label{fig:tucker-trick}
\end{figure}

Relabel the tree $G$ as follows. Replace ${T}_i$ by the core tensor $C_i$ and absorb the matrices $A_i^{(j')}$ into the adjacent vertices. That is, for each adjacent vertex $j'$ connected to $i$ via edge $(i, j')$ with $j' \ne j$, modify the local tensor $ T_{j'}$ by contracting it with the corresponding factor matrix,
\[
 T_{j'} \leftarrow g( T_{j'} \otimes  A_i^{(j')}),
\]
which simply consists of multiplying the matrix $ A_i^{(j')}$ along the appropriate mode of $T_{j'}$. We now have a contraction network of the same structure as initially that represents the same tensor, but one of the dimensions was strictly reduced. This contradicts the assumption that $G$ was minimal. Hence, for a minimal representation, the bond dimension on each edge must match the corresponding component of the multilinear rank of the local tensors.

\textbf{Sufficiency.} Consider the flattening of the full tensor defined by the cut on an arbitrary edge $(a,b)$. After cutting this edge, the tensor $T$ represented by $G$ can be written as $T=g(\widehat T_a\otimes \widehat T_b)$, where $\widehat T_a$ and $\widehat T_b$ are the effective local tensors associated with the subtrees $G_a$ and $G_b$ resulting from the contraction along all internal edges within each subtree, excluding edge $(a,b)$. Without loss of generality, assume that $a$ is the parent of $b$. The flattening matrices of $\widehat T_a$ and $\widehat T_b$ associated with the respective subtree are $\widehat M_a$ and $\widehat M_b$, both of which have rank $r_{ab}$ because of Lemma \ref{lem:edge-rank}. The global flattening matrix $T_{(a,b)}$ resulting from the cut at $(a,b)$ is then
\[
    T_{(a,b)}=\widehat M_a \widehat M_b^\top.
\]
Since $\widehat M_a$ and $\widehat M_b$ have linearly independent columns, $\mathrm{rank}(T_{(a,b)})=r_{ab}$ holds for every edge $(a,b)\in E$. 

Now, suppose that the tuple of bond dimensions $\mathbf r$ fails to be minimal at the edge $(a,b)$. Then, by definition, there exists a smaller tuple $\mathbf s \le \mathbf r$, with at least one strict inequality, such that $ T\in\mathrm{TN}(G, \mathbf s)$. By Lemmas \ref{lem:uniqueTTNrank} and \ref{lem:TTNrank}, there exists an edge $(a', b')$ such that $\mathrm{rank}(T_{(a',b')})\le s_{a'b'}< r_{a'b'}$. This contradicts the rank equalities derived above. Hence, $\mathbf r$ is minimal. 
\end{proof}

\begin{example}
    As an example of (i) how admissibility inequalities are checked, and (ii) how they imply minimality, we revisit the trees of Figure \ref{fig:min_vs_nonmin_tree}.
    
    (i) Let (1) be labelled with bond dimensions $(r_{01}, r_{02}, r_{13}, r_{14})=(6,\, 4,\, 2,\, 3)$. Beginning from the leaves, the admissibility condition at vertex $v_2$ requires $r_{02}\le \dim \mathbb V_{v_2}=2$, so $r_{02}$ is not admissible. Checking $r_{01}$ from $v_1$ gives $r_{01}\le \dim\mathbb V_{v_1} r_{13} r_{14}=1\cdot 2\cdot 3=6$, so $r_{01}$ is admissible from this side. However, checking it from $v_0$ yields the bound $r_{01}\le \dim\mathbb V_{v_0}r_{02}=4$. Therefore the bond dimension tuple is not admissible.
    
    Now, for tree (2), let $(r_{01}, r_{02}, r_{13}, r_{14})=(2, 2, 2, 3)$. Repeating the check as described for (1) confirms that all bond dimensions are admissible.
    
    (ii) Place a tensor with full effective multilinear rank at every vertex in tree (2). Then, for every edge $(i,j)$, the flattening of the local tensor $T_i$ across that edge has rank exactly $r_{ij}$. From Lemma \ref{lem:edge-rank}, the unfolding of the global tensor $T$ on this edge has the same rank $r_{ij}$. Reducing $r_{ij}$ would therefore lower this rank, changing the represented tensor.
\end{example}

Having determined that equality of effective multilinear rank and bond dimensions on all local tensors is equivalent to minimality for an individual tensor, we now show that the set of tensors satisfying this property forms a Zariski open and dense subset of the ambient variety.

\begin{theorem}[Genericity]\label{theorem:genericity}
Let $G=(V,E)$ be a tree and $\mathbf r$ a tuple of bond dimensions. If the tuple $\mathbf r$ is admissible, then $\mathrm{TN}^\circ(G, \mathbf{r})$ is a Zariski open and dense subset of $\mathrm{TN}(G, \mathbf r)$. On the other hand, if the tuple $\mathbf r$ is not admissible, then $\mathrm{TN}^\circ(G, \mathbf r)=\emptyset$.
\end{theorem}

\begin{proof}[Proof]
Assume that $\mathbf{r}$ is admissible. To verify the existence of a tensor in $\mathrm{TN}^\circ(G, \mathbf{r})$, we need to place at each vertex $i \in V$ a tensor with full effective multilinear rank. By Proposition \ref{prop:fullEffectiveRank}, this condition defines a nonempty Zariski open subset of each local space. This proves that $\mathrm{TN}^\circ(G, \mathbf{r}) \neq \emptyset$. By Corollary \ref{cor:TNvariety}, the set $\mathrm{TN}(G, \mathbf{r})$ is an irreducible algebraic variety defined by minors of the edge flattenings. Since $\mathrm{TN}^\circ(G, \mathbf{r})$ is the complement of the finite union of all $\mathrm{TN}(G, \mathbf{s})$ with $\mathbf{s} \le \mathbf{r}$ and $\mathbf r \ne \mathbf s$, and each $\mathrm{TN}(G, \mathbf{s})$ is a variety, $\mathrm{TN}^\circ(G, \mathbf{r})$ is Zariski open. Combined with the existence argument, $\mathrm{TN}^\circ(G, \mathbf{r})$ is a nonempty Zariski open subset of the irreducible variety $\mathrm{TN}(G, \mathbf{r})$, and therefore the former is dense in the latter \cite{idealsvarieties}. 

On the other hand, assume that $\mathbf r$ is not admissible. Then, there exists an edge $(i, j)\in E$ for which the inequality \eqref{eq:admissible} is violated. But then a local tensor $T_i$ can never have $\mathrm{rank}(T_i^{(j)})=r_{ij}$, which by Theorem \ref{theorem:minimality} implies that $\mathbf{r}$ is not minimal for any $T \in \mathrm{TN}(G, \mathbf{r})$, i.e., $\mathrm{TN}^\circ(G, \mathbf{r}) = \emptyset$.
\end{proof}

\begin{remark}
	The second claim also follows from
\cite[Lemma 4.2]{bernardiDimensionTensorNetwork2023}.
\end{remark}

This characterisation provides answers to our main questions for general tree tensor networks, as established in Theorem \ref{theorem:genericity}:

\begin{enumerate}[label=A\arabic*]
    \item \textbf{Minimality:} $\mathrm{TN}^\circ(G, \mathbf{r}) \neq \emptyset$ if and only if $\mathbf{r}$ is admissible.

    \item \textbf{Genericity:} When $\mathrm{TN}^\circ(G, \mathbf{r}) \neq \emptyset$, it forms a nonempty Zariski open and dense subset of $\mathrm{TN}(G, \mathbf{r})$, by Theorem \ref{theorem:genericity}.
\end{enumerate}

\begin{remark}
If the underlying field is $\mathbb{R}$ or $\mathbb{C}$, then $\mathrm{TN}^\circ(G,\mathbf{r})$ is also open and dense in $\mathrm{TN}(G,\mathbf{r})$ with respect to the Euclidean topology. Openness in the Euclidean topology follows immediately from Zariski openness. To establish density, let
\[
T = g(T_1 \otimes \cdots \otimes T_d) \in \mathrm{TN}(G,\mathbf{r}),
\]
where some of the local tensors $T_i$ fail to have full effective multilinear rank. For each $i$, choose a sequence $(T_i^{(t)})_{t \in \mathbb{N}}$ of tensors of full effective multilinear rank such that $T_i^{(t)} \to T_i$ as $t \to \infty$. Define 
\(
T^{(t)} := g\bigl(T_1^{(t)} \otimes \cdots \otimes T_d^{(t)}\bigr).
\)
Then each $T^{(t)} \in \mathrm{TN}^\circ(G,\mathbf{r})$ by Theorem \ref{theorem:minimality}, and by continuity of $g$, we have $T^{(t)} \to T$. Thus $\mathrm{TN}^\circ(G,\mathbf{r})$ is Euclidean open and dense in $\mathrm{TN}(G,\mathbf{r})$.
\end{remark}

\section{Practical reduction of bond dimensions to tree tensor network rank}\label{sec:practical}
Knowing whether a tuple of bond dimensions is minimal is relevant when performing \emph{rank truncation} to obtain a data-sparse tree tensor network approximation. Given a tensor $T$, a tree graph $G=(V,E)$, and bond dimensions $\mathbf{r}$, these methods will approximate $T$ by a tree tensor network $\widetilde{T} \in \mathrm{TN}(G, \mathbf{r})$. However, if the bond dimensions $\mathbf{r}$ are not minimal, then $\widetilde{T}$, contrary to what a user of these methods might be expecting, will have a tree tensor network \emph{rank} strictly smaller than $\mathbf{r}$ because $\mathrm{TN}^\circ(G, \mathbf{r}) = \emptyset$. This implies that a more data-efficient representation of $\widetilde{T}$ is possible with the same approximation error. 

Importantly, the above inconvenience is not always visible for current state-of-the-art algorithms for approximation by tree tensor networks. The standard algorithm, \cite[Algorithm 2]{grasedyck}, for example, mathematically computes a series of truncated HOSVDs \cite{DeLathauwerEtAl} of reshaped tensors starting from the leaves in the tree graph $G$ and moving to the root. Each of these T-HOSVDs will compute several flattenings, compute their singular value decompositions, and then project the flattened tensor onto the dominant $r$-dimensional left singular subspace, where $r$ is the appropriate component from the tuple of bond dimensions. One might hope that if a non-minimal tuple of bond dimensions is provided to the truncation algorithm, it will nevertheless show up somewhere as an unexpected additional number of singular values equal to zero when computing the singular value decomposition of a flattening. However, this is \emph{not} always the case. Consider the next example.

\begin{figure}[tb]
    \centering
    \includegraphics[width=0.75\linewidth]{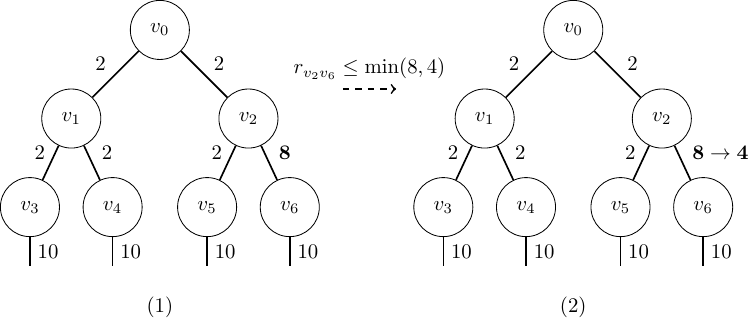}
     \caption{Illustration of non-minimality missed by leaves-to-roots compression. (1) A binary tree representing a tensor $T\in\mathbb R^{10\times10\times10\times10}$ with bond dimension $r_{v_2v_6}=8$. Although the leaf tensor at $v_6$ has full rank towards its parent, the tensor at $v_2$ has effective rank 4 towards this edge, so admissibility requires $r_{v_2v_6}\le4$. (2) Applying \eqref{eq:reduced_r} reduces $8\rightarrow4$, yielding a minimal representation.}
    \label{fig:example_2}
\end{figure}

\begin{example}
Let $T$ be a generic element of $\mathbb{R}^{10 \times 10 \times 10 \times 10}$ and consider the tree graph and bond dimensions in Figure \ref{fig:example_2}(1). Then, all flattenings of $T$ will have the maximal possible rank $\min\{m,n\}$ if the flattening is an $m \times n$ matrix. In the first step, \cite[Algorithm 2]{grasedyck} will truncate a generic $T$ to a tensor with multilinear rank $\mathbf{s}=(2,2,2,8)$ by computing all the $10 \times 1000$ flattening matrices, computing their singular value decompositions, all of whose singular values are strictly positive, and then setting the smallest $10-s_i$ singular values to zero. In the second step, the $2 \times 2 \times 2 \times 8$ core tensor is reshaped to a $4 \times 16$ matrix, which generically has rank equal to $4$. Hence, \cite[Algorithm 2]{grasedyck} will truncate this matrix to bond dimensions $(2,2)$ using a singular value decomposition (which will reveal all $4$ singular values are strictly positive). As both components in this tuple are equal, this will generically correspond to the rank. Finally, \cite[Algorithm 2]{grasedyck} will compute the transfer tensors from the leaves layer to the root layer, though no singular values are computed here. In the end, \cite[Algorithm 2]{grasedyck} will return 4 factor matrices with size $10 \times s_i$, $i=1,\dots,4$, two transfer tensors of size $2 \times 2 \times 2$ for the left branch and $2 \times 8 \times 2$ for the right branch, and a vector for representing the $2\times2$ diagonal matrix at the root. 

At no point in the algorithm was there any indication \emph{in the computed singular values} that the bond dimensions in Figure \ref{fig:example_2}(1) are not minimal. The number of nonzero singular values was always greater than the desired truncation rank from the given bond dimensions. Yet, Theorem \ref{theorem:minimality} shows that the labeling in Figure \ref{fig:example_2}(1) is not minimal. We can, in fact, see this from the $2 \times 8 \times 2$ transfer tensor that is returned by the algorithm: by Carlini and Kleppe's result, the multilinear rank of this tensor will be at most $(2,4,2)$, indicating the bond dimension from the rightmost leaf to its parent can be reduced from 8 to 4 while representing the same tensor. This reduction is shown in Figure \ref{fig:example_2}(2), and we check that the bond dimensions are admissible in this case. By Theorem \ref{theorem:genericity}, we now expect that $T$'s rank truncation has tree tensor network rank equal to the bond dimensions of Figure \ref{fig:example_2}(2). 

\end{example}

This example highlights that bond redundancy can be invisible for state-of-the-art algorithms. However, local conditions from Definition \ref{def:admissible} identify non-minimal bond dimensions directly from the network representation. These conditions can be used to reduce a non-minimal tuple to a minimal tuple as follows. Whenever an admissibility inequality at a vertex $i$ fails for an edge $(i,j)$, replace
\begin{equation}\label{eq:reduced_r}
    r_{ij}\leftarrow \dim\mathbb V_i\prod_{k\in\mathrm{nb}(i)\setminus\{j\}} r_{ik}.
\end{equation}
Repeat this step until condition \eqref{eq:admissible} is satisfied for all edges. 

We show that the above procedure terminates and produces a unique minimal tuple of bond dimensions, independent of the order of reductions. First, we need three lemmas.

\begin{lemma}[Local uniqueness of violation and reduction]\label{lem:local_uniqueness} 
Let $G=(V, E)$ be a tree graph and $\mathbf{r}$ be bond dimensions. Let $i\in V$.
At most one admissibility inequality at $i$ can fail. If such a violation occurs for some edge $(i,j)$, then performing the substitution from \eqref{eq:reduced_r} restores all admissibility inequalities at $i$.
\end{lemma}

\begin{proof}
Let $C=\dim\mathbb V_i$. Suppose admissibility fails at $i$ for some edge $(i,j)$, so that
\[
r_{ij} > C\prod_{k\ne j} r_{ik}.
\]
Since $C\ge 1$ and $r_{ik}\ge 1$, the right-hand side is at least $r_{ij'}$ for every $j'\ne j$, hence $r_{ij} > r_{ij'}$ for all $j'\ne j$. Thus, at most one admissibility inequality can fail. Now define $r'_{ij}=C\prod_{k\ne j} r_{ik}$. Then $r'_{ij}\ge r_{ij'}$ for all $j'\ne j$, so all admissibility inequalities for $j'\ne j$ are automatically satisfied. The inequality for $(i,j)$ holds with equality by construction. Hence all admissibility inequalities at $i$ hold after the replacement.
\end{proof}

\begin{lemma}[Global termination]
    Any sequence of reductions \eqref{eq:reduced_r} terminates after finitely many steps.
\end{lemma}

\begin{proof}
By Lemma \ref{lem:local_uniqueness}, each step reduces exactly one bond dimension of $\mathbf r$. Since $\mathbf r^{(0)}\in \mathbb N^{|E|}$, each $r_{ij}$ can be reduced at most $\Delta_{ij} := r_{ij}^{(0)} - 1$ times. Hence the process terminates, and the total number of reductions is upper bounded by
$
    \sum_{(i,j)\in E} \Delta_{ij}.
$
\end{proof}

\begin{lemma}[Closure under coordinatewise maximum]\label{lem:max_tuple}
Let $\mathbf s,\mathbf t \le \mathbf r^{(0)}$ be admissible bond dimension tuples. Define $\mathbf u$ by
\begin{equation}\label{eq:u_ij}
u_{ij} := \max\{s_{ij},t_{ij}\}
\quad \forall(i,j)\in E.
\end{equation}
Then $\mathbf u$ is also admissible.
\end{lemma}

\begin{proof}
Let $C_i=\dim\mathbb V_i$. Fix a vertex $i$ and a neighbour $j\in\mathrm{nb}(i)$.
Since $\mathbf s$ and $\mathbf t$ are admissible, we have
\[
s_{ij} \le C_i \prod_{k\ne j} s_{ik},
\quad
t_{ij} \le C_i \prod_{k\ne j} t_{ik}.
\]
By definition,
$ u_{ij} = \max\{s_{ij},t_{ij}\}.$
In particular, for each $k\ne j$,
$ u_{ik} = \max\{s_{ik},t_{ik}\} \ge s_{ik} $ and $u_{ik} \ge t_{ik}$.
Hence,
\[
    \prod_{k\ne j} u_{ik}
    \ge
    \prod_{k\ne j} s_{ik},
    \quad
    \prod_{k\ne j} u_{ik}
    \ge
    \prod_{k\ne j} t_{ik}.
\]
Multiplying by $C_i$ yields
\[
    C_i\prod_{k\ne j} u_{ik}
    \ge
    C_i\prod_{k\ne j} s_{ik}\ge s_{ij},
    \quad
    C_i\prod_{k\ne j} u_{ik}
    \ge
    C_i\prod_{k\ne j} t_{ik}\ge t_{ij}.
\]
Therefore,
\[
u_{ij}
=
\max\{s_{ij},t_{ij}\}
\le
C_i\prod_{k\ne j} u_{ik}.
\]
Since this holds for all $i$ and $j$, the tuple $\mathbf u$ is admissible.
\end{proof}

With these lemmas, we can now prove correctness.

\begin{theorem}[Global uniqueness of reduction terminus]
    For any initial bond dimension tuple $\mathbf{r}^{(0)}$, all reduction sequences terminate at the same admissible tuple.
\end{theorem}

\begin{proof}
    Let
    \[
    S=\{\mathbf v \text{ admissible} : \mathbf v \le \mathbf r^{(0)}\}.
    \]
    By Lemma \ref{lem:max_tuple}, the finite set $S$ is closed under coordinatewise maxima, and hence has a unique maximal element. Denote this element by $\mathbf u$. We claim that every tuple $\mathbf r^{(t)}$ produced by the reduction algorithm satisfies $\mathbf r^{(t)} \ge \mathbf u$. This is shown by induction on $t$. For $t=0$ the claim holds because $\mathbf u \le \mathbf r^{(0)}$ by definition. Suppose $\mathbf r^{(t)} \ge \mathbf u$. If $r^{(t)}_{ij}$ is reduced, then
    \[
    r^{(t+1)}_{ij}
    = C\prod_{k\ne j} r^{(t)}_{ik}
    \ge C\prod_{k\ne j} u_{ik}
    \ge u_{ij},
    \]
    where the last inequality follows from admissibility of $\mathbf u$. All other components are unchanged, so $\mathbf r^{(t+1)} \ge \mathbf u$. Thus $\mathbf r^{(t)} \ge \mathbf u$ for all $t$. When the algorithm terminates we obtain an admissible tuple $\mathbf r^{(t)}$. Hence $\mathbf r^{(t)} \in S$, while also $\mathbf r^{(t)} \ge \mathbf u$. By maximality of $\mathbf u$ in $S$, we conclude $\mathbf r^{(t)}=\mathbf u$. Therefore every reduction sequence terminates at the same admissible tuple $\mathbf u$.
\end{proof}

In particular, if $\mathrm{TN}^\circ(G, \mathbf{r}^{(0)})=\emptyset$, then every tensor in $\mathrm{TN}(G, \mathbf{r}^{(0)})$ has tree tensor network rank strictly smaller than $\mathbf{r}^{(0)}$. Thus truncation algorithms operating within a non-minimal model may not reveal the redundancy in the bond dimensions. Minimality should therefore be verified prior to truncation, which can be done using the intrinsic reduction procedure established above.

\section*{Acknowledgments}
J.~J.~and N.~V.~were supported by the FWO Weave Project ``Group-invariant neural tensor train networks'' with number G011624N. T.~S.~was supported by the FWO Junior Postdoctoral Fellowship with number 1219723N. We thank the two reviewers for their careful reading and constructive comments, which have helped improve the clarity and scope of this manuscript. In particular, we have expanded Section \ref{sec:practical} to better develop the practical implications of our results, and clarified the relation to existing literature in Section \ref{sec:lit_review}. 

\printbibliography

@article{carlinikleppe,
    title = {Ranks derived from multilinear maps},
    journal = {Journal of Pure and Applied Algebra},
    volume = {215},
    number = {8},
    pages = {1999-2004},
    year = {2011},
    doi = {10.1016/j.jpaa.2010.11.010},
    author = {Enrico Carlini and Johannes Kleppe},}

@article{grasedyck,
author = {Grasedyck, Lars},
title = {Hierarchical Singular Value Decomposition of Tensors},
journal = {SIAM Journal on Matrix Analysis and Applications},
volume = {31},
number = {4},
pages = {2029-2054},
year = {2010},
doi = {10.1137/090764189}
}

@book{idealsvarieties,
  title     = {Ideals, Varieties, and Algorithms},
  subtitle  = {An Introduction to Computational Algebraic Geometry and Commutative Algebra},
  author    = {Cox, David A. and Little, John and O'Shea, Donal},
  series    = {Undergraduate Texts in Mathematics},
  edition   = {4},
  publisher = {Springer Cham},
  year      = {2015},
  doi       = {10.1007/978-3-319-16721-3},
  isbn      = {978-3-319-16720-6},
  pages     = {XVI, 646},
  issn      = {0172-6056},
  eissn     = {2197-5604}
}

@book{hackbusch2019tensor,
  title={Tensor Spaces and Numerical Tensor Calculus},
  author={Hackbusch, Wolfgang},
  year={2019},
  publisher={Springer}
}

@misc{limandye,
      title={Tensor network ranks}, 
      author={Ke Ye and Lek-Heng Lim},
      year={2019},
      note={arXiv:2301.12345}
}

@article{bernardiDimensionTensorNetwork2023,
  title = {Dimension of {{Tensor Network}} Varieties},
  author = {Bernardi, Alessandra and Lazzari, Claudia De and Gesmundo, Fulvio},
  date = {2023-12},
  journal = {Communications in Contemporary Mathematics},
  volume = {25},
  number = {10},
  pages = {2250059},
  doi = {10.1142/S0219199722500596}
}

@article{landsbergGeometryTensorNetwork2012,
  author    = {Joseph M. Landsberg and Yang Qi and Ke Ye},
  title     = {On the geometry of tensor network states},
  journal   = {Quantum Information \& Computation},
  volume    = {12},
  number    = {3-4},
  pages     = {346--354},
  year      = {2012},
  publisher = {Rinton Press},
  doi       = {10.5555/2230976.2230988},
}

@article{buczynskaHackbuschConjectureTensor2015,
  title = {The {{Hackbusch}} Conjecture on Tensor Formats},
  author = {Weronika Buczy{\'n}ska and Jaros{\l}aw Buczy{\'n}ski and Mateusz Micha{\l}ek},
  date = {2015-10},
  journal = {Journal de Mathématiques Pures et Appliquées},
  volume = {104},
  number = {4},
  pages = {749--761},
  doi = {10.1016/j.matpur.2015.05.002},
}

@article{Sidiropoulos_2017,
   title={Tensor Decomposition for Signal Processing and Machine Learning},
   volume={65},
   DOI={10.1109/tsp.2017.2690524},
   number={13},
   journal={IEEE Transactions on Signal Processing},
   publisher={Institute of Electrical and Electronics Engineers (IEEE)},
   author={Sidiropoulos, Nicholas D. and De Lathauwer, Lieven and Fu, Xiao and Huang, Kejun and Papalexakis, Evangelos E. and Faloutsos, Christos},
   year={2017},
   month=jul, pages={3551–3582} }

@inproceedings{stoundenmire,
 author = {Stoudenmire, Edwin and Schwab, David J},
 booktitle = {Advances in Neural Information Processing Systems},
 editor = {D. Lee and M. Sugiyama and U. Luxburg and I. Guyon and R. Garnett},
 publisher = {Curran Associates, Inc.},
 title = {Supervised Learning with Tensor Networks},
 volume = {29},
 year = {2016}
}

@book{Landsberg2012,
  author    = {J. M. Landsberg},
  title     = {Tensors: Geometry and Applications},
  year      = {2012},
  series    = {Graduate Studies in Mathematics},
  volume    = {128},
  publisher = {American Mathematical Society},
  address   = {Providence, RI},
  isbn      = {978-0-8218-6907-9},
  pages     = {xx+439}
}

@book{Greub1978,
  author       = {Werner Greub},
  title        = {Multilinear Algebra},
  year         = {1978},
  edition      = {2},
  series       = {Universitext},
  volume       = {136},
  publisher    = {Springer New York},
  address      = {New York, NY},
  isbn         = {978-0-387-90284-5},
  doi          = {10.1007/978-1-4613-9425-9},
  pages        = {VIII, 296}
}

@article{ShiDuanVidal2006,
  title = {Classical simulation of quantum many-body systems with a tree tensor network},
  author = {Shi, Y.-Y. and Duan, L.-M. and Vidal, G.},
  journal = {Phys. Rev. A},
  volume = {74},
  issue = {2},
  pages = {022320},
  numpages = {4},
  year = {2006},
  month = {Aug},
  publisher = {American Physical Society},
  doi = {10.1103/PhysRevA.74.022320}
}

@article{Verstraete01032008,
author = {F. Verstraete and V. Murg and J.I. Cirac},
title = {Matrix product states, projected entangled pair states, and variational renormalization group methods for quantum spin systems},
journal = {Advances in Physics},
volume = {57},
number = {2},
pages = {143--224},
year = {2008},
publisher = {Taylor \& Francis},
doi = {10.1080/14789940801912366}
}

@article{evenblyTensorNetworkStates2011,
  title = {Tensor {{Network States}} and {{Geometry}}},
  author = {Evenbly, G. and Vidal, G.},
  year = 2011,
  month = nov,
  journal = {Journal of Statistical Physics},
  volume = {145},
  number = {4},
  pages = {891--918},
  issn = {1572-9613},
  doi = {10.1007/s10955-011-0237-4}
}

@article{bachmayrTensorNetworksHierarchical2016a,
  title = {Tensor {{Networks}} and {{Hierarchical Tensors}} for the {{Solution}} of {{High-Dimensional Partial Differential Equations}}},
  author = {Bachmayr, Markus and Schneider, Reinhold and Uschmajew, André},
  date = {2016-12-01},
  journal = {Foundations of Computational Mathematics},
  volume = {16},
  number = {6},
  pages = {1423--1472},
  issn = {1615-3383},
  doi = {10.1007/s10208-016-9317-9},}

@article{Haegeman2014GeometryMPS,
  author  = {Haegeman, Jutho and Mariën, Michaël and Osborne, Tobias J. and Verstraete, Frank},
  title   = {Geometry of matrix product states: Metric, parallel transport, and curvature},
  journal = {Journal of Mathematical Physics},
  volume  = {55},
  number  = {2},
  pages   = {021902},
  year    = {2014},
  doi     = {10.1063/1.4862851}
}

@article{holtzManifoldsTensorsFixed2012,
  title = {On Manifolds of Tensors of Fixed {{TT-rank}}},
  author = {Holtz, Sebastian and Rohwedder, Thorsten and Schneider, Reinhold},
  year = 2012,
  month = apr,
  journal = {Numerische Mathematik},
  volume = {120},
  number = {4},
  pages = {701--731},
  issn = {0029-599X, 0945-3245},
  doi = {10.1007/s00211-011-0419-7}
}

@article{USCHMAJEW2013133,
title = {The geometry of algorithms using hierarchical tensors},
journal = {Linear Algebra and its Applications},
volume = {439},
number = {1},
pages = {133-166},
year = {2013},
issn = {0024-3795},
doi = {10.1016/j.laa.2013.03.016},
author = {André Uschmajew and Bart Vandereycken},}

@article{Orus2014,
	author = {Or{\'u}s, Rom{\'a}n},
	date = {2014-11-26},
	doi = {10.1140/epjb/e2014-50502-9},
	id = {Or{\'u}s2014},
	isbn = {1434-6036},
	journal = {The European Physical Journal B},
	number = {11},
	pages = {280},
	title = {Advances on tensor network theory: symmetries, fermions, entanglement, and holography},
	volume = {87},
	year = {2014}}

@article{grasedyckIntroductionHierarchicalRank2011,
  title = {An {{Introduction}} to {{Hierarchical}} ({{H-}}) {{Rank}} and {{TT-Rank}} of {{Tensors}} with {{Examples}}},
  author = {Grasedyck, Lars and Hackbusch, Wolfgang},
  year = 2011,
  journal = {Computational Methods in Applied Mathematics},
  volume = {11},
  number = {3},
  pages = {291--304},
  issn = {1609-9389, 1609-4840},
  doi = {10.2478/cmam-2011-0016},
}

@article{hackbuschNewSchemeTensor2009a,
  title = {A {{New Scheme}} for the {{Tensor Representation}}},
  author = {Hackbusch, W. and K{\"u}hn, S.},
  year = 2009,
  month = oct,
  journal = {Journal of Fourier Analysis and Applications},
  volume = {15},
  number = {5},
  pages = {706--722},
  issn = {1531-5851},
  doi = {10.1007/s00041-009-9094-9},}

@article{SCHOLLWOCK201196,
  title = {The Density-Matrix Renormalization Group in the Age of Matrix Product States},
  author = {Schollw{\"o}ck, Ulrich},
  year = 2011,
  journal = {Annals of Physics},
  volume = {326},
  number = {1},
  pages = {96--192},
  issn = {0003-4916},
  doi = {10.1016/j.aop.2010.09.012},}

@article{novikov2022,
  title = {Automatic Differentiation for {R}iemannian Optimization on Low-Rank Matrix and Tensor-Train Manifolds},
  author = {Novikov, Alexander and Rakhuba, Maxim and Oseledets, Ivan},
  year = 2022,
  journal = {SIAM Journal on Scientific Computing},
  volume = {44},
  number = {2},
  pages = {A843-A869},
  doi = {10.1137/20M1356774},}

@incollection{Grasedyck2024,
  title = {Stable Truncation and Root-Independent Normalization of Tree Tensor Networks},
  booktitle = {Multiscale, Nonlinear and Adaptive Approximation {{II}}},
  author = {Grasedyck, Lars and Kr{\"a}mer, Sebastian and Moser, Dieter},
  editor = {DeVore, Ronald and Kunoth, Angela},
  year = 2024,
  pages = {231--265},
  publisher = {Springer Nature Switzerland},
  address = {Cham},
  doi = {10.1007/978-3-031-75802-7_12},
  isbn = {978-3-031-75802-7}
}

@mastersthesis{adams2024,
  author  = {Jonas Adams},
  title   = {Group-Invariant Tree Tensor Networks},
  school  = {KU Leuven},
  year    = {2024},
  type    = {Master's thesis},
  note     = {\nopunct\url{https://repository.teneo.libis.be/delivery/DeliveryManagerServlet?dps\_pid=IE19320203}}
  }

@article{DeLathauwerEtAl,
    author = {De Lathauwer, Lieven and De Moor, Bart and Vandewalle, Joos},
    title = {On the Best Rank-1 and Rank-(R1 ,R2 ,. . .,RN) Approximation of Higher-Order Tensors},
    journal = {SIAM Journal on Matrix Analysis and Applications},
    volume = {21},
    number = {4},
    pages = {1324-1342},
    year = {2000},
    doi = {10.1137/S0895479898346995}
}
\end{document}